\documentclass{article}
\usepackage{amsmath}
\newcommand{\ebinom}[2]
{\binom{#1}{#2}}
\begin{document}
\title{On the expansion of the power of any polynomial
$(1+x+x^2+x^3+x^4+\textrm{ etc.})^n$\footnote{Delivered to the
St.--Petersburg Academy July 6, 1778. Originally published as
{\em De evolutione potestatis polynomialis
cuiuscunque $(1+x+x^2+x^3+x^4+\textrm{ etc.})^n$},
Nova Acta Academiae Scientarum Imperialis Petropolitinae \textbf{12} (1801),
47-57, and reprinted in
\emph{Leonhard Euler, Opera Omnia}, Series 1:
Opera mathematica, Volume 16, Birkh\"auser, 1992.
A copy of the original text is available
electronically at the Euler Archive, at http://www.eulerarchive.org. This paper
is E709 in the Enestr\"om index.}}
\author{Leonhard Euler\footnote{Date of translation: May 20, 2004.
Translated from the Latin by Jordan
Bell, 3rd year undergraduate in Honours Mathematics, School of Mathematics
and Statistics, Carleton University, Ottawa, Ontario, Canada. Email:
jbell3@connect.carleton.ca.
This translation was written
during an NSERC USRA supervised by Dr. B. Stevens.}}
\date{}

\maketitle

1. We begin with a power of the binomial $(1+x)^n$, for which,
having been expanded, by usual
custom we designate the coefficients
of each power $x^\lambda$ by the character $\ebinom{n}{\lambda}$,
\footnote{Translator: Euler in fact uses the symbol
$\left( \begin{array}{r}
\alpha\\
\hline
\beta
\end{array} \right)$ for $\ebinom{\alpha}{\beta}$} so that
it is thus:
\[
(1+x)^n=1+\ebinom{n}{1}x +\ebinom{n}{2}xx+\ebinom{n}{3}x^3+\ebinom{n}{4}x^4+
\textrm{ etc.} \ldots \ebinom{n}{n}x^n,
\]
for which it will then be:
\begin{align*}
\ebinom{n}{1}=n; \quad \ebinom{n}{2}=\frac{n(n-1)}{1\cdot 2}; \quad
\ebinom{n}{3}=\frac{n(n-1)(n-2)}{1\cdot 2\cdot 3}; \\
\ebinom{n}{4}=\frac{n(n-1)(n-2)(n-3)}{1\cdot 2\cdot 3\cdot 4} \quad
\textrm{ etc.},
\end{align*}
and in general:
\[
\ebinom{n}{\lambda}=\frac{n(n-1)(n-2)\cdots (n-\lambda+1)}{1\cdot 2\cdot 3
\cdots \lambda};
\]
from this it follows for the case $\lambda=0$ and $\lambda=n$ to
be $\ebinom{n}{0}=\ebinom{n}{n}=1$, and truly also in general
$\ebinom{n}{\lambda}=\ebinom{n}{n-\lambda}$. 
Besides this in fact, it willa be helpful to note that in the cases
in which $\lambda$ is a negative number or in which it is a number larger
than $n$, it is designated for the form $\ebinom{n}{\lambda}$ always
to be equal to nothing.

2. Since this reckoning is assisted and shortened to no small degree
by means of these
characters, we will also use similar characters in the expansion of
trionomial powers, quadrinomials, and generally for any polynomial whatsoever.
Henceforth, to this end, for the above
characters, we will adjoin to those used for binomials the exponent 2,
seeing that in this situation no ambiguity is to be feared, since in these calculations
it is customary for no powers to occur;
in this manner, for the expansion of binomial powers we will have:
\[
(1+x)^2=\ebinom{n}{0}^2+\ebinom{n}{1}^2x+\ebinom{n}{2}^2xx+\ebinom{n}{3}^2x^3+
\textrm{ etc.},
\]
for which indeed it should be remembered to be in general
$\ebinom{n}{\lambda}^2=\ebinom{n}{n-\lambda}^2$,
and then also always $\ebinom{n}{0}^2=\ebinom{n}{2}^2=1$,
and as well for those forms to disappear into nothing 
in which $\lambda$ is either an integral number that is negative or
positive and greater than $n$.

3. These same characters will therefore be used for the expansion
of any polynomial whatsoever, providing that for trinomials we adjoin
exponents of threes, for quadrinomials of fours, for quintinomials
of fives, and so forth, namely in this way:
 
For trinomials $(1+x+xx)^n$, this expansion unfolds:
\[
\ebinom{n}{0}^3+\ebinom{n}{1}^3x+\ebinom{n}{2}^3x^2+\ebinom{n}{3}^3x^3+
\ebinom{n}{4}^3x^4+\textrm{ etc.}
\]

For quadrinomials $(1+x+xx+x^3)^n$, this expansion unfolds:
\[
\ebinom{n}{0}^4+\ebinom{n}{1}^4x+\ebinom{n}{2}^4x^2+\ebinom{n}{3}^4x^3+
\ebinom{n}{4}^4x^4+\textrm{ etc.}
\]

For quintinomials $(1+x+xx+x^3+x^4)^n$, this expansion unfolds:
\[
\ebinom{n}{0}^5+\ebinom{n}{1}^5x+\ebinom{n}{2}^5xx+\ebinom{n}{3}^5x^3+
\ebinom{n}{4}^5x^4+\ebinom{n}{5}^5x^5+\textrm{ etc.}
\]

4. We will now inquire into unraveling the actual values of the
characters inscribed with the exponents $3, 4, 5, 6$, etc., and
we will see in what way these characters will be able to be
determined by those inscribed with 2, for which of course the
signification is most well known. Therefore
we run through each of these cases of powers of polynomials in order.

\begin{center} 
{\Large The expansion of the power of the trionomial $(1+x+xx)^n$}
\end{center} 

5. We represent the series that arises in this way: 
\[
\ebinom{n}{0}^3+\ebinom{n}{1}^3x+\ebinom{n}{2}^3xx+\ebinom{n}{3}^3x^3+
\ebinom{n}{4}^3x^4+\textrm{ etc.},
\]
of which the last term will be equal to $\ebinom{n}{2n}^3x^{2n}$,
for which as we already knew, the coefficient
$\ebinom{n}{2n}^3$ is to be equal to unity, and so too
the first term $\ebinom{n}{0}^3$. Next indeed, because the coefficients
from the end procede as the same as those in order, it then follows
to be:
\[
\ebinom{n}{1}^3=\ebinom{n}{2n-1}^3; \ebinom{n}{2}^3=\ebinom{n}{2n-2}^3,
\]
in fact, in general $\ebinom{n}{\lambda}^3=\ebinom{n}{2n-\lambda}^3$.
Again it is evident for the form $\ebinom{n}{\lambda}^3$ to take
the value of nothing in the cases in which $\lambda$ is an integral
number that is negative and in those cases in which it is positive and larger
than $2n$.

6. However, before we take up the determination of these characters,
the expansion of simpler cases before our eyes
will be by no means unfitting:

\begin{tabular}{r|l}
$n$&$(1+x+xx)^n$\\
\hline
0&1\\
1&$1+x+xx$\\
2&$1+2x+3xx+2x^3+x^4$\\
3&$1+3x+6xx+7x^3+6x^4+3x^5+x^6$\\
4&$1+4x+10xx+16x^3+19x^4+16x^5+10x^6+4x^7+x^8$\\
5&$1+5x+15xx+30x^3+45x^4+51x^5+45x^4+30x^3+$ etc. [sic]\\
6&$1+6x+21xx+50x^3+90x^4+126x^5+141x^6+126x^7+$ etc.\\
&etc.
\end{tabular}

From the last case, in which $n=6$, it is therefore apparent to be:
\begin{align*}
\ebinom{6}{0}^3=1; \ebinom{6}{1}^3=6; \ebinom{6}{2}^3=21; \ebinom{6}{3}^3=50;
\ebinom{6}{4}^3=90;\\
\ebinom{6}{5}^3=126; \ebinom{6}{6}^3=141; \ebinom{6}{7}^3=126;
\ebinom{6}{8}^3=90;\\
\ebinom{6}{9}^3=30; \ebinom{6}{10}^3=21; \ebinom{6}{11}^3=6; \ebinom{6}{12}^3=1.
\end{align*}

7. So that we can now investigate the way in which the characters
that arise from the trionomials may be derived according to the similar
characters from the binomials from before; we will represent the given
powers under the form of binomials in this way: $[1+x(1+x)]^n$,
for which the expansion therefore produces this progression:
\begin{align*}
1+\ebinom{n}{1}^2x(1+x)+\ebinom{n}{2}^2xx(1+x)^2+\ebinom{n}{3}^2x^3(1+x)^3+\\
\ebinom{n}{4}^2x^4(1+x)^4+
\textrm{ etc.},
\end{align*}
in which the general term will have this form: $\ebinom{n}{\alpha}x^\alpha
(n-x)^\alpha$. 

8. We now consider the expansion of any given power $x^\lambda$ of $x$,
of which the coefficient is $\ebinom{n}{\lambda}^3$, whose value we will
investigate. To this end, the power $x^\lambda$ ought to be produced
from each one of the binomial members made in this way,
so that it will be formed by them. On the other hand, the general
form is $\ebinom{n}{\alpha}x^\alpha (1+x)^\alpha$, from which:
\[
(1+x)^\alpha =1+\ebinom{\alpha}{1}^2x+\ebinom{\alpha}{2}^2x^2+
\ebinom{\alpha}{3}^2x^3+\ebinom{\alpha}{4}^2x^4+\textrm{ etc.}
\]
Because in this place
 generally the term that occurs is $\ebinom{\alpha}{\beta}^2x^\beta$,
this having been extended with $\ebinom{n}{\alpha}^2x^\alpha$ 
produces $\ebinom{\alpha}{\beta}^2\ebinom{n}{\alpha}^2x^{\alpha+\beta}$.
For if therefore it will have been $\alpha+\beta=\lambda$,
the coefficient $\ebinom{\alpha}{\beta}^2\ebinom{n}{\alpha}^2$
will be equal to the coefficient sought from $\ebinom{n}{\lambda}^3$.

9. In order to elicit the value of the coefficient $\ebinom{n}{\lambda}^3$,
it is necessary precisely for the letters $\alpha$ and $\beta$
to take all integral value which are able to provide $\alpha+\beta=\lambda$.
Moreover, it is evident for both of these numbers $\alpha$ and $\beta$ to
be taken as either negative or as larger than $n$,
in which cases this form vanishes; then also, if it is $\beta>\alpha$,
the form $\ebinom{\alpha}{\beta}^2$ will be equal to nothing.
Therefore then, the maximum value that will be assumed for $\alpha$ is
$\lambda$, at which time indeed $\beta$; from this it follows:
\begin{align*}
\left. \begin{array}{rc|c|c|c|c|c}
\textrm{if}&\alpha&\lambda&\lambda-1&\lambda-2&\lambda-3&\lambda-4\\ 
\textrm{it will be}&\beta&0&1&2&3&4
\end{array} \right\} \textrm{ etc.}
\end{align*}

10. Therefore it is attained that for each case
that arises, the value of the character $\ebinom{n}{\lambda}^3$
that is being searched for will be given collected together in one sum,
such that this determination is born:
\begin{align*}
\ebinom{n}{\lambda}^3=\ebinom{\lambda}{0}^2\ebinom{n}{\lambda}^2+
\ebinom{\lambda-1}{1}^2\ebinom{n}{\lambda-1}^2+
\ebinom{\lambda-2}{2}^2
\ebinom{n}{\lambda-2}^2+\\
\ebinom{\lambda-3}{3}^2\ebinom{n}{\lambda-3}^2+
\textrm{ etc.},
\end{align*}
and in this way this value is expressed by known parts, of which the
number for any case is finite.

11. So that this will be well understood, we expand these simple
cases in which $\lambda$ itself takes the values $0, 1, 2, 3, 4$, etc.,
and it will be as follows:
\begin{tabular}{l}
$\ebinom{n}{0}^3=1; \ebinom{n}{1}^3=\ebinom{1}{0}^2\ebinom{n}{1}=n$;\\
$\ebinom{n}{2}^3=\ebinom{2}{0}^2\ebinom{n}{2}^2+\ebinom{1}{1}^2\ebinom{n}{1}^2=
\ebinom{n}{2}^2+\ebinom{n}{1}^2=\frac{n(n-1)}{1\cdot 2}+n=
\frac{n(n+1)}{1\cdot 2}$;\\
$\ebinom{n}{3}^3=\ebinom{3}{0}^2\ebinom{n}{3}^2+\ebinom{2}{1}^2\ebinom{n}{2}^2$,
that is,\\
$\ebinom{n}{3}^3=\ebinom{n}{3}^2+2\ebinom{n}{2}^2$,\\
$\ebinom{n}{4}^3=\ebinom{4}{0}^2\ebinom{n}{4}^2+\ebinom{3}{1}^2\ebinom{n}{3}^2+
\ebinom{2}{2}^2\ebinom{n}{2}^2$, that is,\\
$\ebinom{n}{4}^3=\ebinom{n}{4}^2+3\ebinom{n}{3}^2+\ebinom{n}{2}^2$,\\
$\ebinom{n}{5}^3=\ebinom{5}{0}^2\ebinom{n}{5}^2+\ebinom{4}{1}^2\ebinom{n}{4}^2+
\ebinom{3}{2}^2\ebinom{n}{3}^2$, that is,\\
$\ebinom{n}{5}^3=\ebinom{n}{5}^2+4\ebinom{n}{4}^2+3\ebinom{n}{3}^2$,\\
$\ebinom{n}{6}^3=\ebinom{6}{0}^2\ebinom{n}{6}^2+\ebinom{5}{1}^2\ebinom{n}{5}^2+
\ebinom{4}{2}^2\ebinom{n}{4}^2+\ebinom{3}{3}^2\ebinom{n}{3}^2$, that is,\\
$\ebinom{n}{6}^3=\ebinom{n}{6}^2+5\ebinom{n}{5}^2+6\ebinom{n}{4}^2+
\ebinom{n}{3}^2$,\\
$\ebinom{n}{7}^3=\ebinom{7}{0}^2\ebinom{n}{7}^2+\ebinom{6}{1}^2\ebinom{n}{6}^2+
\ebinom{5}{2}^2\ebinom{n}{5}^2+\ebinom{4}{3}^2\ebinom{n}{4}^2$, that is,\\
$\ebinom{n}{7}^3=\ebinom{n}{7}^2+6\ebinom{n}{6}^2+10\ebinom{n}{5}^2+
4\ebinom{n}{2}^2$,\\
$\ebinom{n}{8}^3=\ebinom{8}{0}^2\ebinom{n}{8}^2+\ebinom{7}{1}^2\ebinom{n}{7}^2+
\ebinom{6}{2}^2\ebinom{n}{6}^2+\ebinom{5}{3}^2\ebinom{n}{5}^2+
\ebinom{4}{4}^2\ebinom{n}{4}^2$, that is,\\
$\ebinom{n}{8}^3=\ebinom{n}{8}^2+7\ebinom{n}{7}^2+15\ebinom{n}{6}^2+
10\ebinom{n}{5}^2+\ebinom{n}{4}^2$,\\
$\ebinom{n}{9}^3=\ebinom{9}{0}^2\ebinom{n}{9}^2+\ebinom{8}{1}^2\ebinom{n}{8}^2+
\ebinom{7}{2}^2\ebinom{n}{7}^2+\ebinom{6}{3}^2\ebinom{n}{6}^2+
\ebinom{5}{4}^2
\ebinom{n}{5}^2$, that is,\\
$\ebinom{n}{9}^3=\ebinom{n}{9}^2+8\ebinom{n}{8}^2+21\ebinom{n}{7}^2+
20\ebinom{n}{6}^2+5\ebinom{n}{5}^2$,\\
$\ebinom{n}{10}^3=\ebinom{10}{0}^2\ebinom{n}{10}^2+\ebinom{9}{1}^2\ebinom{n}{9}^2+
\ebinom{8}{2}^2\ebinom{n}{8}^2+\ebinom{7}{3}^2\ebinom{n}{7}^2+
\ebinom{6}{4}^2\ebinom{n}{6}^2+\ebinom{5}{5}^2\ebinom{n}{5}^2$,
that is,\\
$\ebinom{n}{10}^3=\ebinom{n}{10}^2+9\ebinom{n}{9}^2+28\ebinom{n}{8}^2+
35\ebinom{n}{7}^2+15\ebinom{n}{6}^2+\ebinom{n}{5}^2$,\\
etc.
\end{tabular}

12. We will now employ this for the purpose of an example for the case
$n=6$, which of course we have expanded above in \S 6, and we will discover:

\begin{tabular}{l}
$\ebinom{6}{0}^3=1$,\\
$\ebinom{6}{1}^3=6$,\\
$\ebinom{6}{2}^3=21$,\\
$\ebinom{6}{3}^3=\ebinom{6}{3}^2+2\ebinom{6}{2}^2=50$,\\
$\ebinom{6}{4}^3=\ebinom{6}{4}^2+3\ebinom{6}{3}^2+\ebinom{6}{2}^2=
15+3\cdot 20+15=90$,\\
$\ebinom{6}{5}^3=\ebinom{6}{5}^2+4\ebinom{6}{4}^2+3\ebinom{6}{3}^2=6+4\cdot 15+
3\cdot 20=126$,\\
$\ebinom{6}{6}^3=\ebinom{6}{6}^2+5\ebinom{6}{5}^2+6\ebinom{6}{4}^2+
\ebinom{6}{3}^2=1+5\cdot 6+6\cdot 15+20=141$,\\
$\ebinom{6}{7}^3=6\ebinom{6}{6}^2+10\ebinom{6}{5}^2+4\ebinom{6}{4}^2+
\ebinom{6}{3}^2$, that is,\\
$\ebinom{6}{7}^3=6+10\cdot 6+4\cdot 15=126$,
\end{tabular}

One may know, by means of it being $\ebinom{6}{\alpha}^3=\ebinom{6}{12-\alpha}^3$,
that it will indeed be $\binom{6}{7}^3=\ebinom{6}{5}^3=126$;
in a similar way, it will be $\ebinom{6}{8}^3=\ebinom{6}{4}^3=90$;
$\ebinom{6}{9}^3=\ebinom{6}{3}^3=50$;
$\ebinom{6}{10}^3=\ebinom{6}{2}^3=21$; $\ebinom{6}{11}^3=\ebinom{6}{1}^3=6$;
and finally $\ebinom{6}{12}^3=\ebinom{6}{0}^3=1$, 
whose values agree completely with those given before.

\begin{center}
{\Large The expansion of the power of the quadrinomial
$(1+x+xx+x^3)^n$}
\end{center}

13. We will then represent the value having been expanded thus:
\[
1+\ebinom{n}{1}^4x+\ebinom{n}{2}^4xx+\ebinom{n}{3}^4x^3+\ebinom{n}{4}^4x^4+
\ebinom{n}{5}^4x^5+\textrm{ etc.,}
\]
where it is of course $\ebinom{n}{0}^4=1$. Next, because the last term is
$x^3n$ it will be $\ebinom{n}{3n}^4=1$, and because the coefficients written
from the end follow as the very same as those in order, it will be
$\ebinom{n}{3n-1}^4=\ebinom{n}{1}^4$, and indeed in general
$\ebinom{n}{3n-\lambda}^4=\ebinom{n}{\lambda}^4$. In this it is
observed that in the cases in which $\lambda$ is an integral
number that is negative, or positive and larger than $3n$,
the values of these forms are to disappear into nothing.
For expressing these, it has been proposed by me to seek out the way in which
these quaternary characters may be expressed by either binary or ternary
characters which have been expressed, in as much as these are known,
 and which can be so determined.

14. Before we take up this task, we put simple cases of the given
form joined together in this table before our eyes:

\begin{tabular}{l|l}
$n$&$(1+x+xx+x^3)^n$\\
\hline
0&1\\
1&$1+x+xx+x^3$\\
2&$1+2x+3xx+4x^3+3x^4+2x^5+x^6$\\
3&$1+3x+6xx+10x^3+12x^4+12x^5+10x^6+6x^7+3x^8+x^9$\\
4&$1+4x+10xx+20x^3+31x^4+40x^5+44x^6+40x^7+31x^8+$ etc.\\
5&$1+5x+15xx+35x^3+65x^4+101x^5+135x^6+155x^7+155x^8+$ etc.\\
6&$1+6x+21xx+56x^3+120x^4+216x^5+$ etc.\\
&etc.
\end{tabular}

15. We now refer to the given form under this binomial:
$[1+x(1+x+xx)]^n$, whose expansion by us presents this series:
\[
1+\ebinom{n}{1}^2[x(1+x+xx)]+\ebinom{n}{2}^2x^2(1+x+xx)^2+
\textrm{ etc.},
\]
whose general term is $\ebinom{n}{\alpha}^2x^\alpha (1+x+xx)^\alpha$.
Now indeed, because $(1+x+xx)^\alpha$ is a power of a trinomial, it will be:
\[
(1+x+xx)^\alpha=1+\ebinom{\alpha}{1}^3 x + \ebinom{\alpha}{2}^3 xx+
\ebinom{\alpha}{3}^3 x^3+ \textrm{ etc.,}
\]
where for another time the general term is $\ebinom{\alpha}{\beta}^3 x^\beta$,
from which it follows for the powers $x^\lambda$ to come into being
as $\lambda=\alpha+\beta$, from which this member
$\ebinom{\alpha}{\beta}^3\ebinom{n}{\alpha}^2 x^\lambda$ is seen
for this power.

16. Therefore in the expansion that is inquired into, the coefficient
of the power $x^\lambda$ will be $\ebinom{n}{\lambda}^4$,
whose value will be revealed if, since it is $\lambda=\alpha+\beta$,
all values of the form $\ebinom{n}{\alpha}^2 \ebinom{\alpha}{\beta}^3$
are gathered together into one sum, which will be formed thus:
\begin{align*}
\ebinom{n}{\lambda}^4=\ebinom{n}{\lambda}^2 \ebinom{\lambda}{0}^3+
\ebinom{n}{\lambda-1}^2 \ebinom{\lambda-1}{1}^3 +
\ebinom{n}{\lambda-2}^2 \ebinom{\lambda-2}{2}^3 +\\
\ebinom{n}{\lambda-3}^2 \ebinom{\lambda-3}{3}^3 
\textrm{ etc.}
\end{align*}
In such a way, all the quaternary characters which are expressed
may be determined
by either the binary or ternary which have been expressed, and which
are already known;
therefore, it is clearly evident that we may write
the numbers $0,1,2,3,4$, etc. in place of $\lambda$, and we wll
obtain:

\begin{tabular}{l}
$\ebinom{n}{0}^4=\ebinom{n}{0}^2 \ebinom{0}{0}^3=1$,\\
$\ebinom{n}{1}^4=\ebinom{n}{1}^2 \ebinom{1}{0}^3=n$,\\
$\ebinom{n}{2}^4=\ebinom{n}{2}^2 \ebinom{2}{0}^3+\ebinom{n}{1}^2 \ebinom{1}{1}^3 = \ebinom{n}{2}^2+n$,\\
$\ebinom{n}{3}^4=\ebinom{n}{3}^2 \ebinom{3}{0}^3 + \ebinom{n}{2}^2 \ebinom{2}{1}^3 + \ebinom{n}{1}^2 \ebinom{1}{2}^3$, that is\\
$\ebinom{n}{3}^4=\ebinom{n}{3}^2+2\ebinom{n}{2}^2+\ebinom{n}{1}^2$,\\
$\ebinom{n}{4}^4=\ebinom{n}{4}^2 \ebinom{4}{0}^3 +
\ebinom{n}{3}^2 \ebinom{3}{1}^3 + \ebinom{n}{2}^2 \ebinom{2}{2}^3 +\ebinom{n}{1}^2 \ebinom{1}{3}^3$, that
is\\
$\ebinom{n}{4}^4=\ebinom{n}{4}^2+ 3\ebinom{n}{3}^2+3\ebinom{n}{2}^2$,\\
$\ebinom{n}{5}^4=\ebinom{n}{5}^2 \ebinom{5}{0}^3 + \ebinom{n}{4}^2 \ebinom{4}{1}^3 + \ebinom{n}{3}^2 \ebinom{3}{2}^3 + \ebinom{n}{2}^2 \ebinom{2}{3}^3$,\\
$\ebinom{n}{6}^4=\ebinom{n}{6}^2 \ebinom{6}{0}^3 + \ebinom{n}{5}^2 \ebinom{5}{1}^3 + \ebinom{n}{4}^2 \ebinom{4}{2}^3 + \ebinom{n}{3}^2 \ebinom{3}{3}^3 +
\ebinom{n}{2}^2 \ebinom{2}{4}^3+$ etc.\\
$\ebinom{n}{7}^4=\ebinom{n}{7}^2 \ebinom{7}{0}^3 + \ebinom{n}{6}^2 \ebinom{6}{1}^3 + \ebinom{n}{5}^2 \ebinom{5}{2}^3 + \ebinom{n}{4}^2 \ebinom{4}{3}^3 +
\ebinom{n}{3}^2 \ebinom{3}{4}^3+$ etc.\\
$\ebinom{n}{8}^4=\ebinom{n}{8}^2 \ebinom{8}{0}^3 +
\ebinom{n}{7}^2 \ebinom{7}{1}^3 + \ebinom{n}{6}^2 \ebinom{6}{2}^3 +
\ebinom{n}{5}^2 \ebinom{5}{3}^3 +\ebinom{n}{4}^2 \ebinom{4}{4}^3 + \ebinom{n}{3}^2 \ebinom{3}{5}^3+$ etc.\\
etc.
\end{tabular}

\begin{center}
{\Large The expansion of the power of the quintinomial
$(1+x+xx+x^3+x^4)^n$}
\end{center} 

17. Thus we present the value of this expanded as so:
\[
1+\ebinom{n}{1}^5 x+\ebinom{n}{2}^5 x^2+ \ebinom{n}{3}^5 x^3+ \ebinom{n}{4}^5 x^4+
\ebinom{n}{5}^5 x^5 + \textrm{ etc.},
\]
 where it is $\ebinom{n}{0}^5=\ebinom{n}{4n}^5=1$, and moreover in general
$\ebinom{n}{\lambda}^5=\ebinom{n}{4n-\lambda}^5$;
then indeed it follows for those values to vanish in the cases in which
$\lambda$ is an integral number which is negative or in which it is positive
and larger than $4n$.

18. Now this very form will be represented as the
binomial $[1+x(1+x+xx+x^3)]^n$, for the expansion of which, in general
a member comes forth as $\ebinom{n}{\alpha}^2 x^\alpha (1+x+xx+x^3)^\alpha$,
where the factor $(1+x+xx+x^3)^\alpha$ includes the term
$\ebinom{\alpha}{\beta}^4 x^\beta$,
which having been joined together, the term
$\ebinom{n}{\alpha}^2 \ebinom{\alpha}{\beta}^4 x^{\alpha+\beta}$
is obtained. Hence if it were $\alpha+\beta=\lambda$,
for a member the coefficient of the power $x^\lambda$ will
be $\ebinom{n}{\alpha}^2 \ebinom{\alpha}{\beta}^4$. At this point,
for the letters $\alpha$ and $\beta$ are taken all the values which they are
able to receive, starting from $\alpha=\lambda$, and so
the coefficient that is being searched for will be:
\begin{align*}
\ebinom{n}{\lambda}^5=\ebinom{n}{\lambda}^2 \ebinom{\lambda}{0}^4+
\ebinom{n}{\lambda-1}^2 \ebinom{\lambda-1}{1}^4 + 
\ebinom{n}{\lambda-2}^2
\ebinom{\lambda-2}{2}^4 +\\
 \ebinom{n}{\lambda-3}^2 \ebinom{\lambda-3}{3}^4+\textrm{ etc.},
\end{align*}
and so all characters expressed by the number 5 will be defined by
those preceding characters expressed by the number 4, along with those
characters expressed by the number 2.

\begin{center}
{\Large General conclusion}
\end{center} 

19. From this, it is now quite clear that if a power of a polynomial
is given in general, comprised of by the number $\theta+1$ of terms, namely
$(1+x+xx+x^3 \ldots x^\theta)^n$, then for those terms with the power
$x^\lambda$, the adjoined coefficient will be
$\ebinom{n}{\lambda}^{\theta+1}$, thus expressed from $\theta$ number of
characters, will be comprised so that it is:
\[
\ebinom{n}{\lambda}^{\theta+1}=\ebinom{n}{\lambda}^2 \ebinom{\lambda}{0}^\theta+
\ebinom{n}{\lambda-1}^2 \ebinom{\lambda-1}{1}^\theta +
\ebinom{n}{\lambda-2}^2 \ebinom{\lambda-2}{2}^\theta+\textrm{ etc.,}
\]
whose form includes all the preceding in it. For example, if we begin
with the value $\theta=1$, the case will be realized of the power
of the binomial $(1+x)^n$, and moreover the characters expressed by unity
are seen from the power of the monomial $1^n$, from which comes
$\ebinom{n}{0}=1$, where indeed all the others depart into nothing. Then,
for the purpose of these cases proceding, we will have that it follows:

\begin{tabular}{l}
$\ebinom{n}{\lambda}^2=\ebinom{n}{\lambda}^2 \ebinom{\lambda}{0}=
\ebinom{n}{\lambda}^2$,\\
$\ebinom{n}{\lambda}^3=\ebinom{n}{\lambda}^2 \ebinom{\lambda}{0}^2+
\ebinom{n}{\lambda-1}^2 \ebinom{\lambda-1}{1}^2 + \ebinom{n}{\lambda-2}^2
\ebinom{\lambda-2}{2}^2 +$ etc.\\
$\ebinom{n}{\lambda}^4=\ebinom{n}{\lambda}^2 \ebinom{\lambda}{0}^3 +
\ebinom{n}{\lambda-1}^2 \ebinom{\lambda-1}{1}^3 + \ebinom{n}{\lambda-2}^2
\ebinom{\lambda-2}{2}^3+$ etc.\\
$\ebinom{n}{\lambda}^5=\ebinom{n}{\lambda}^2 \ebinom{\lambda}{0}^4+
\ebinom{n}{\lambda-1}^2 \ebinom{\lambda-1}{1}^4 + \ebinom{n}{\lambda-2}^2 
\ebinom{\lambda-2}{2}^4+$ etc.\\
$\ebinom{n}{\lambda}^6=\ebinom{n}{\lambda}^2 \ebinom{\lambda}{0}^5 +
\ebinom{n}{\lambda-1}^2 \ebinom{\lambda-1}{1}^5 + \ebinom{n}{\lambda-2}^2
\ebinom{\lambda-2}{2}^5+$ etc.\\
$\ebinom{n}{\lambda}^7=\ebinom{n}{\lambda}^2 \ebinom{\lambda}{0}^6 + 
\ebinom{n}{\lambda-1}^2 \ebinom{\lambda-1}{1}^6 + \ebinom{n}{\lambda-2}^2
\ebinom{\lambda-2}{2}^6+$ etc.\\
etc.
\end{tabular}

\end{document}